\documentclass[12pt,fleqn]{article}

\usepackage{amsmath}
\usepackage{amsfonts}
\usepackage{amsthm}
\newcommand{\Rset}{\mathbb{R}}
\newcommand{\Cset}{\mathbb{C}}

\begin{document}

\title{GENERALIZED L\'EVY STOCHASTIC AREAS AND SELFDECOMPOSABILITY}
\author{Zbigniew J. Jurek\footnote{This work was completed while
the author was visiting (2002/03) Wayne State University, Detroit,
USA.} , Wroc\l aw, Poland.}
\date{\emph{Stat. \& Probab. Letters}\textbf{64}(2003),213-222.}

\maketitle

\newtheorem{thm}{THEOREM}
\newtheorem{lem}{LEMMA}
\newtheorem{prop}{PROPOSITION}
\newtheorem{cor}{COROLLARY}

\theoremstyle{remark}
\newtheorem{rem}{REMARK}

\begin{quote}
\noindent {\footnotesize \textbf{ABSTRACT.} We show that a
conditional characteristic function of generalized L\'evy
stochastic areas can be viewed as a product a selfdecomposable
distribution (i.e., L\'evy class L distribution) and its
background driving
characteristic function.
This provides a stochastic interpretation for a ratio of some
Bessel functions as well as examples of characteristic functions
from van Dantzig class.}

\end{quote}

\noindent \textbf{1. An introduction.} For a planar Brownian
motion $\textbf{B}_t = ( Z_t,\Tilde{Z}_t)$ and a stochastic area
process
\[
\mathcal{A}_u = \int_{0}^{u} Z_sd\Tilde{Z}_s - \Tilde{Z}_sdZ_s ,
\qquad u>0,
\]
Paul L\'evy (1950) proved that its conditional characteristic
function is of the form
\begin{equation*}
\mathbb{E}[e^{it \mathcal{A}_u} | \textbf{B}_u=a] =
\frac{tu}{\sinh tu} \exp[- \frac{|a|^2}{2u}(tu \coth tu-1)],\ \
t\in \Rset,
\end{equation*}
where $a \in \Rset^2$ and $u \ge 0$ are fixed. Thus, in
particular,
\begin{equation}
\mathbb{E}[e^{it \mathcal{A}_u} | \textbf{B}_u=(\sqrt{u},\sqrt{u}
)] = \frac{tu}{\sinh tu} \exp[-(tu \coth tu -1)],\ \ t\in \Rset.
\end{equation}
For another proof of L\'evy stochastic area formula and related topics
involving more stochastic argument we refer to
Williams (1976) and Yor (1980).

In Jurek (2001), p.248, it was showed that in the formula (1) one
has a product of the \emph{class $L$ characteristic function}
$bt/\sinh bt$, ($b \in \Rset$ is a fixed parameter), and  a
characteristic function $\exp[-(bt\coth bt -1)]$ of its
\emph{background driving random variable} that is from the class
$\mathcal{U}$ of s-selfdecomposable distributions. Consequently,
by Proposition 3, in Iksanov, Jurek and Schreiber (2003), we have
that (1) is again in class $L$. In this note we prove that a
similar result holds for a generalized L\'evy stochastic areas;
cf.Theorem 1. Here hyperbolic trigonometric functions are
replaced by the modified Bessel function of first kind $I_{\nu}$. As a
by-product we get a formula for a ratio of such functions and
provide examples of van Dantzig class $\mathcal{D}$; cf.
Proposition 1, Corollaries 2-4 and the formula (12).

Bessel functions $I_{\nu}(t)$, as completely monotone functions in
$\sqrt{t}$ or $\sqrt{\nu}$ or as functions of a pair $(t,\nu)$,
appeared in many areas of probability and can be interpreted by
some stochastic functionals; cf. Ismail and Kelker (1979), but in
particular, cf. Pitman and Yor (1981) and references there-in. In
this note we provide stochastic way of looking at $I_{\nu}(t)$,
and their ratios, as function in $t\in \Rset$.

Our proofs here rely on different representations of Bessel
functions and therefore for ease of reading, in the Appendix,
we have collected all needed formulas. This also makes this paper
more self-contained.
\medskip

\textbf{2. Selfdecomposable and s-selfdecomposable distributions.}

\emph{Selfdecomposable} random variables $X$ (or their probability
distribution or their characteristic functions $\phi$) are defined
as  limits of normalized partial sums of infinitesimal independent
random variables. Thus they are infinitely divisible, in short:
$X\in ID$. Let  $L$ stand for the class of all either
selfdecomposable random variables $X$ or their probability
distributions $\mu$ or their characteristic functions $\phi$. Then
\begin{equation}
\phi\in L\quad \mbox{iff}\quad \forall{( 0<c<1)}\exists
(char.f.\,\, \psi_{c}) \  \ \phi(\cdot) =
\phi(c\cdot)\psi_{c}(\cdot),
\end{equation}
and this is so called \emph{the selfdecomposability property.} One
of the major characterizations of the class $L$, called \emph{the random
integral representation} is as follows
\begin{equation}
X\in L\quad \mbox{iff}\quad X \stackrel{d}{=} \int_{(0,\infty)}
e^{-s}dY(s),\ \mbox{with}\ E[\log(1+ |Y(1)|)] <\infty,
\end{equation}
where  $\stackrel{d}{=}$ means equality in distribution and
$Y(\cdot)$ is a L\'evy process called \underline{b}ackground
\underline{d}riving \underline{L}\'evy \underline{p}rocess; in
short: $Y$ is BDLP of $X$ or $X$ is \emph{driven} by $Y$; cf.
Jurek and Mason (1993), Chapter 3. For BDLP identification
purposes note that if $\psi$ is the characteristic function of
$Y(1)$ in (3) then
\begin{equation}
\psi(t) = \exp[t(\log \phi(t))'], \ \ \ t\neq 0, \ \ \mbox{and} \
\ \psi(0)=1.
\end{equation}
In particular, class $L$ characteristic functions are differentiable for $t\ne 0$
and $\psi$ is infinitely divisible characteristic function with a
finite logarithmic moment; cf. Jurek (2001), Proposition 3.
Finally let us recall that $(L,\ast)$ is a closed topological
subsemigroup of the semigroup $ID$.
\newline
Class $L$ is quite rich and contains among others: all stable laws,
gamma $\gamma_{\alpha,\lambda}$, Laplace (double exponential),
$t$-Student, $F$-Fisher, log-normal, log-F, $\log|t|$, inverse
Gaussian, Barndorff-Nielsen generalized hyperbolic, generalized
gamma, etc.

Another class of infinitely divisible distributions needed here is the class
$\mathcal{U}$ of \emph{s-selfdecomposable probability measures}.
One of the descriptions of $\mathcal{U}$ is as follows
\begin{equation}
\psi\in \mathcal{U}\,\,\, \mbox{iff}\,\,\, \forall{( 0<c<1)}\exists
\mbox{(char.f.}\,\,\rho_{c}\in ID)\,\, \psi(t) =
\psi^c(c\,t)\rho_{c}(t), \,\, t \in \Rset,
\end{equation}
where the power is well defined as $\psi(t)\neq 0$. To the above we refer to as
\emph{the s-selfdecomposability probability.} From (5), in particular,
we infer that $(\mathcal{U},\ast)$ is also a closed convolution
topological semigroup. In fact, we have the inclusions $ L\subset
\mathcal{U}\subset ID$. Similarly, for s-selfdecomposable
distributions we have analogous random integral representation
\begin{equation}
X \in \mathcal{U}\ \ \mbox{iff}\ \ X \stackrel{d}{=}\int_{(0,1)}s\, dY(s),
\end{equation}
and as before we refer to $Y$ as the BDLP for s-selfdecomposable $X$. From (6) we
infer that
if $\phi$ and $\psi$ are characteristic functions of $X$ and $Y(1)$, respectively, then
\[
\psi(t)=\exp[(t \log\phi(t))'], \ \ t \neq 0 \ \ \mbox{and} \ \
\psi(0)=1.
\]

For those and other relations between the classes $L$ and
$\mathcal{U}$ we refer to Jurek (1985) or Iksanov, Jurek and
Schreiber (2003).
\medskip

\textbf{3. Stochastic areas and selfdecomposability. } For
$p>-1/2$, let us define a Gaussian process $\textbf{V}^p$ and the
corresponding  generalized L\'evy stochastic area as follows
\[
\textbf{V}^p_t=(V_t^p,\Tilde{V}_t^p):=t^{-p}\int_0^t
s^p\,d\textbf{B}_s,\ \ \ \ \mbox{and}\ \ \mathcal{A}^p:=\int_0^1
V_s^p\, d\Tilde{V}_s^p-\Tilde{V}_s^p\,dV_s^p,
\]
where $\textbf{B}_s,\, s\ge 0,$ is the planar Brownian motion. In
Biane and Yor (1987) (see also Yor (1989), p. 3052 or Duplantier
(1989)) the following generalization of L\'evy stochastic area
formula was proved
\begin{equation}
\mathbb{E}[e^{i\lambda
\mathcal{A}^p}|\textbf{V}^p_1=a]=\frac{|\lambda|^{\nu}}{2^{\nu}\Gamma(\nu+1)I_{\nu}(|\lambda|)}
\exp[-\frac{|a|^2}{2}\,|\lambda|\,
\frac{I_{\nu+1}(|\lambda|)}{I_{\nu}(|\lambda|)}],
\end{equation}
where $I_{\nu}(z)$ is the modified Bessel function of the index
$\nu:=p +1/2>0$ and $a\in \Rset^2$ is fixed. In the 'products' in
(7) we recognize a similar formula as in (1) and , moreover, the
components have analogous interpretations. Namely we have the
following result.
\begin{thm}
(i) For $\nu >-1$, functions
\[
B_{\nu}(t):= \frac{t^{\nu}}{2^{\nu}\Gamma(\nu+1)\,I_{\nu}(t)}, \ \
t\in\Rset,\ \ \mbox{are class L characteristic functions.}
\]
More explicitly, these are the characteristic functions of random
variables of the form
\[
X_{\nu}:=\sum_{k=1}^{\infty}\,\,z_{\nu,k}^{-1}\,\, \eta_k, \ \
\mbox{where}\ \ \eta_k \ \ \mbox{are i. i. d. Laplace rv's,}
\]
and $z_{\nu,k}$, $k=1,2,...$, are the zeros of the function
$z^{-\nu}J_{\nu}(z)$.  Here $J_{\nu}$ is  the Bessel function of the
first kind. Furthermore
\[
\log B_{\nu}(t)= - \int_{\Rset\setminus\{0\}}(1-\cos tx)
(\sum_{k=1}^{\infty}\,e^{-|z_{\nu,k}||x|})/|x|dx.
\]

(ii) Functions
\[
b_{\nu}(t):= \exp\big[-t\,\frac{I_{\nu+1}(t)}{I_{\nu}(t)}\big], \
\ t\in \Rset,\ \ \mbox{are class $\mathcal{U}$ characteristic
functions.}
\]
Moreover, these are the background driving characteristic
functions for $B_{\nu}$, i.e., $\log b_{\nu}(t)= t\, d(\log
B_{\nu}(t))/dt, \ \ t\neq 0$. Furthermore,
\begin{equation}
(t/2)\,I_{\nu+1}(t)= I_{\nu}(t)\,\int_0^{\infty}(1-\cos
tx)(\sum_{k=1}^{\infty}\,|z_{\nu,k}|e^{-|z_{\nu,k}|x})dx.
\end{equation}

 (iii) For all positive constants $\alpha$ and $\beta$,
 the functions $B_{\nu}^{\alpha}(t)\cdot b_{\nu}^{\beta}(t)$ are
s-selfdecomposable (class $\mathcal{U}$ distributions). Moreover,
for $0<\beta \le \alpha$ the resulting characteristic functions are
selfdecomposable (class $L$ characteristic functions).

(iv) For $\nu>-1/2$, functions
$1/B_{\nu}(it)=\frac{\Gamma(\nu+1)J_{\nu}(t)}{(t/2)^{\nu}},
\ \ t \in \Rset$, are characteristic
functions corresponding to the probability densities
\[
f_{\nu}(x):=
\frac{\Gamma(\nu+1)}{\Gamma(\frac{1}{2})\Gamma(\nu+\frac{1}{2})}\,(1-x^2)^{\nu-\frac{1}{2}}1_{[-1,1]}(x).
\]
These are not infinitely divisible characteristic functions.
\end{thm}

\emph{Proof.} From A.3(a) in the Appendix, we have
\begin{equation}
I_{\nu}(z)=\frac{(z/2)^{\nu}}{\Gamma(\nu+1)}
\prod_k^{\infty}\big(1+ \frac{z^2}{z_{\nu,k}^2} \big), \ \ \ z \in
\Cset,\ \ \ \nu>-1,
\end{equation}
where $z_{\nu,k}$ are real zeros of the function
$z^{-\nu}J_{\nu}(z)$, ordered accordingly
to their absolute values; cf. also A.1(a) in the Appendix.
Since Laplace (also called the double
exponential) random variable $\eta$ has characteristic function
\[
\phi_{\eta}(t)=
\int_{-\infty}^{\infty}e^{itx}2^{-1}e^{-|x|}dx=(1+t^2)^{-1}=
\exp\int_{\Rset\setminus\{0\}}(\cos\,tx-1)\frac{e^{-|x|}}{|x|}dx,
\]
we easily conclude that it is class L distribution. Furthermore,
\begin{equation}
B_{\nu}(t)=\prod_{k=1}^{\infty}\big(1+\frac{t^2}{z_{\nu,k}^2}\big)^{-1}
=\exp\int_{\Rset\setminus\{0\}}(\cos tx
-1)(\sum_{k=1}^{\infty}\,e^{-|z_{\nu,k}||x|})/|x|\,dx
\end{equation}
cf. Jurek (1996), p.177 for more details. Since class $L$ is
closed under dilations (multiplication of rv by scalars) and the
weak convergence therefore $X_{\nu}$ has class $L$ distribution as
well, with the characteristic function $B_{\nu}(t)$,which proves
the part (i).

Part (ii). By  Proposition 3 in Jurek (2001) or see (4), for any class $L$
characteristic function $\phi(t)$ we have that $\psi(t):=exp(t
d(log \phi(t))/dt)$ exists (for $t\neq 0$) and it is the
characteristic function of BDLP. Since $z\,\frac{d}{dz}\,I_{\nu}(z)-\nu I_{\nu}(z)=
z\,I_{\nu+1}(z)$, by A.8, we get
\[
t\,d(\log B_{\nu}(t))/dt = \nu-t\, \frac{I_{\nu}'(t)}{I_{\nu}(t)}= -t\,
\frac{I_{\nu+1}(t)}{I_{\nu}(t)}=\log b_{\nu}(t),
\]
which proves that $b_{\nu}(t)$ are indeed characteristic functions
of the BDLP for rv $X_{\nu}$. Using (4) and Jurek (1996), we get
the formula (8). Finally, since $X_{\nu}$ are given by series on
independent Laplace rv therefore by Proposition 3 in Iksanov,
Jurek and Schreiber (2003) to infer that $b_{\nu}(t)$ is from the
class $\mathcal{U}$ of s-selfdecomposable distributions.

Part (iii). Since $L\subset \mathcal{U}$ and both classes are
closed under taking positive powers therefore
$B_{\nu}^{\alpha}(t)\cdot b_{\nu}^{\beta}(t)\in \mathcal{U}$,
which proves the first part of (iii). Let $\alpha\ge \beta>0$.
Since $b_{\nu}(t)\in \mathcal{U}$ therefore $B_{\nu}(t)\cdot
b_{\nu}(t) \in L$, by Theorem 1 in Iksanov, Jurek and Schreiber
(2003). Since also $B_{\nu}^{\frac{\alpha}{\beta} -1}(t) \in L$ and
$L$ is a semigroup, we conclude that
$B_{\nu}(t)^{\alpha/\beta}\cdot b_{\nu}(t)\in L$. Taking powers we
get the second part of (iii).

For part (iv), using A.5 and part (i) of the theorem, we get
\begin{equation}
\frac{1}{B_{\nu}(it)}=
\frac{1}{B(\frac{1}{2},\nu+\frac{1}{2})}\int_{-1}^1
e^{itx}(1-x^2)^{\nu-\frac{1}{2}}dx,
\end{equation}
where $B(x,y)=\frac{\Gamma(x)\Gamma(y)}{\Gamma(x+y)}$ is Euler's
beta function (Euler's integral of the first kind). Since, by A.1,
$J_{\nu}(t)$ has real zeros therefore it can not be
infinite divisible. This completes the proof of the Theorem.

For some related information on random variables of the form $X_{\nu}$
and the corresponding Dirichlet series we refer to Jurek (2000). Here, in particular, we have
the following conditions on  zeros of Bessel functions.

\begin{cor}
For zeros $z_{\nu,k}$, $k=1,2,...$, of the function
$z^{-\nu}J_{\nu}(z)$ with $\nu>-1$,  we have
\begin{equation*}
\sum_k|z_{\nu,k}|\int_{(|x|>1)}\log |x|\,e^{-|z_{\nu,k}||x|}dx<\infty,
\end{equation*}
\begin{equation*}
\sum_k|z_{\nu,k}|\int_{(|x|>0)}\log (1+ |x|^2)\,e^{-|z_{\nu,k}||x|}dx<\infty.
\end{equation*}
\end{cor}
[Both integrals can not be expressed in simple functions. The first involves the
incomplete gamma function and the second sine and cosine integrals;
cf. Gradshetyn and Ryzhik (1965) formulae 4.358 and 4.338, respectively.]

\emph{Proof.} From (i) in Theorem 1 we infer that $X_{\nu}$ has L\'evy spectral $M$ of the form
\[
M(A)=\sum_{k=1}^{\infty}\int_A e^{-|z_{\nu,k}||x|}\,|x|^{-1}dx, \
\ \mbox{for all Borel subsets of} \ \Rset\setminus\{0\}.
\]
And from random integral representation (3) or from (8) in Theorem 1(ii) we have
\[
M(A)=\int_{(0,\infty)}N(e^sA)ds \ \ \mbox{and} \ \ N(A)=\sum_{k=1}^{\infty}\int_A |z_{\nu,k}|e^{-|z_{\nu,k}||x|}dx,
\]
and $N$ is a L\'evy spectral measure of BDLP $Y$ with finite logarithmic moment on $\{|x|>1\}$.
This gives the first condition in Corollary. For the second one we use the fact that
\[
\int_{\{|x|>0\}}\frac{x^2}{1+x^2}M(dx)<\infty,
\]
and the integration by parts formula.

\begin{rem}
Biane and Yor (1987) expanded Brownian path along Legendre
polynomials and computed conditional characteristic function of
stochastic area given values of few first coefficients. In that
case again we have formulas involving the products (7); cf.
formula (4.10) in Yor (1989).
\end{rem}
\begin{rem}
The fractions $I_{\nu+1}/I_{\nu}$ of Bessel functions appear in
the formula for the background driving characteristic functions
$b_{\nu}$ for the selfdecomposable characteristic functions
$B_{\nu}(t)$; cf. Theorem 1(ii). In Jurek (2001), in Example 2, it
was found that functions
\[
t\to \exp[-|t|\,\frac{K_{\nu-1}(|t|)}{K_{\nu}(|t|)}\,],\ \, \ \ \
t\in \Rset,
\]
involving fractions of modified Bessel functions  $K_{\nu}$, are
background driving characteristic functions as well. This time
they correspond to Student's t-distributions with $2\nu$ degrees
of freedom.
\end{rem}
\medskip
\textbf{4. van Dantzig distributions and selfdecomposability.}

Recall that an analytic characteristic function $\phi(z)$
defined in a strip $-a< \Im z< b$, for some $a,b,>0$, belongs to
\emph{van Danztig class} $\mathcal{D}$, if $\psi(t):=1/\phi(it)$
can be extended to a characteristic function; cf. Lukacs (1968).
Equivalently,
\begin{equation}
\psi(t) \cdot \phi(it)=1, \ \ \ t \in \Rset, \ \ ( \mbox{van Dantzig pair} \ \ (\phi,\psi)).
\end{equation}
If $\psi(t)=\phi(t)$ we say $\phi$ is \emph{self-reciprocal characteristic function}
in  van Danzitg class $\mathcal{D}$. The three elementary
examples of pairs $(\psi,\phi)$ are:
\[
 (\cos t,(\cosh t)^{-1}), \ \ \ \Big(\frac{\sin t}{t}, \frac{t}{\sinh t}\Big)\ \ \mbox{and}\ \ (e^{-ct^2},
e^{-ct^2})\ \ \mbox{with} \ \  c\ge0.
\]
The normal characteristic function is an example of
self-reciprocal one. Here is another example.
\begin{prop}
For $\nu>-1/2$, the ratios
$\frac{{J_{\nu}(t)}}{I_{\nu}(t)}, \ \ t \in \Rset,$ are
non-infinitely-divisible self-reciprocal characteristic functions
from van Danztig class $\mathcal{D}$.
\end{prop}
\emph{Proof.} From  Theorem 1, parts (i) and (iv), we have that
the ratio is indeed a characteristic function. Since it has real
zeros therefore $\frac{{J_{\nu}(t)}}{I_{\nu}(t)}$ can
not be infinitely divisible. Its self-reciprocal property follows
from the identities A.2 in the Appendix and the fact that
$B_{\nu}$ are even functions; cf. in the Appendix formulae A.3 or
A.3a. Thus the proof is complete.

\begin{cor}
All characteristic functions $B_{\nu}$, for $\nu>-1/2$, belong to
van Dantzig class $\mathcal{D}$.
\end{cor}
This is a consequence of Theorem 1, parts (i) and (iv). Here are
more explicit formulae of class L and van Dantzig characteristic
functions.
\begin{cor}
For $n=0,1,2,...$ and $t\in \Rset$, functions
\begin{multline*}
B_{n+\frac{1}{2}}(t)=\frac{((2n+1)!!)^{-1}\,\,t^{2n+1}}{\big(\sum_{\{k:\,n-k\,even\}}\,a_{k,n}t^{n-k}\big)\,\sinh
t + \big(\sum_{\{k:\,n-k\,odd\}}\,a_{k,n}t^{n-k}\big)\,\cosh t}
\end{multline*}
where
$a_{k,n}:=(-1)^k\frac{(n+k)!}{2^k\,k!(n-k)!}$, $0\le k\le n$, are
characteristic functions of class $L$ distributions that belong to
van Dantzig class $\mathcal{D}$.
\end{cor}

The proof is a combination of the formulae
\[
I_{(n+\frac{1}{2})}(z)=\sqrt{\frac{2}{\pi}}z^{-(n+\frac{1}{2})}
\]
\[
\cdot \big[\big(\sum_{\{k:n-k\,even\}}a_{k,n}\,z^{n-k}\big)\,\sinh z +
\big(\sum_{\{k:n-k\,odd\}}a_{k,n}\,z^{n-k}\big)\,\cosh z \big]
\]
from A.6a and
$\Gamma(n+\frac{1}{2})=\frac{(2n-1)!!}{2^n}\sqrt{\pi}$ $, \ n \ge
1,$ and $\Gamma(\frac{1}{2})=\sqrt{\pi}$, with the part (i) of
Theorem 1.
\begin{cor}
For $n=0,1,2,...$ and van Dantzig class $\mathcal{D}$
characteristic functions $B_{n+\frac{1}{2}}(t)$, we have the
equality
\begin{equation*}
B_{n+\frac{1}{2}}(t)=\frac{1} {(2n+1)!!
\,\frac{d^n}{(t\,dt)^n}\big(\frac{\sinh t}{t}\big)}\, \, , t \in
\Rset,
\end{equation*}
where $\frac{d^n}{(t\,dt)^n}$ is n-times composition of the differential operator
$ f \to \frac{1}{t}\,f'(t)$.
\end{cor}
This is a consequence of the identity
$I_{n+\frac{1}{2}}(z)=z^{n+\frac{1}{2}}\sqrt{\frac{2}{\pi}}
\frac{d^n}{(z\,dz)^n}\big(\frac{\sinh z}{z}\big)$, cf. A.7a, in
the Appendix, and the Theorem 1(i).
\begin{rem} Note that our class $L$ characteristic functions
$B_{\frac{1}{2}}(t),B_{\frac{3}{2}}(t)$ and $B_{\frac{5}{2}}(t)$
are exactly the examples (of the van Dantzig class $\mathcal{D}$)
given in Lukacs (1968), on p.121, formulae (14.2a),(14.3a) and
(14.4a). However, they were obtained there by different inductive
procedures.
\end{rem}

\begin{rem}
It might be worthy to notice that Lukacs' 1968 examples of class
$\mathcal{D}$ characteristic functions, given below his formula
(13) on p. 120, are class L characteristic functions. The same
applies to those given at the bottom of page 121 (convolutions of
normal distributions with several gamma distributions). All in all
it might be true that in each non self-reciprocal pair
$(\phi,\psi)$ of characteristic functions satisfying van Dantzig
condition (12) one of them is selfdecomposable. Or if it is not
the case, then characterize that part of $\mathcal{D}$ that possesses
this property.
\end{rem}

\begin{rem}
In search for examples of van Danztig pairs one may follow recent paper by de Meyer,
Roynette, Vallois and Yor (2002). In Theorem 4.1 Authors have the equality
\[
\ \ \  \ \mathbb{E}[\exp\,i\lambda\textbf{B}_T]\,\mathbb{E}[\exp\frac{\lambda^2}{2}T]=1, \ \  \lambda \in \Rset,
\]
for a standard Brownian motion $\textbf{B}_t, t\ge 0$, and some independent
of it stopping times $T$.
\end{rem}

\textbf{Acknowledgement.} Author would like to thank Professor Marc Yor,
from University of Paris VI, for some valuable comments and for providing the references [14],
[16] and [17].

\medskip
\textbf{5. Appendix.} Recall that the Bessel functions are the
solutions of the differential equation
\[
\frac{d^2Z_{\nu}}{dz^2}+\frac{1}{z}\,\frac{dZ_{\nu}}{dz}+(1-\frac{\nu^2}{z^2})Z_{\nu}=0.
\]
Our main interest here is in the Bessel function of the first kind
$J_{\nu}(z)$ and their modified versions $I_{\nu}(z)$,
with $z\in \Cset$.

For ease of reading we collect here some formulas for Bessel
functions. They are either taken directly from Gradshteyn and
Ryzhik (1965), and in that case they are labelled by bold face
numbers, or are combination of those.

\begin{itemize}
\item[A.1]
$J_{\nu}(z)=\frac{(z/2)^{\nu}}{\Gamma(\nu+1)}
\prod_{m=1}^{\infty}\big(1- \frac{z^2}{z_{\nu,m}^2}\big), \ \ \nu
\neq -1,-2,-3,...,$ \ \ \ \ \ \ \ \ \ \ \ \ \ \textbf{8.544}

\mbox{where} \ \ $z_{\nu,m}$ \ \ \mbox{are zeros of the function}
\ \ $z^{-\nu}J_{\nu}(z)$, \ \mbox{ordered } \\
\mbox{by absolute value of their real parts}.

\item[A.1a] \mbox{For} \ $\nu >-1$ \mbox{all zeros} \ $z_{\nu,m}$\
\  \mbox{are real numbers.} \ \ \ \ \ \  \ \ \ \ \ \ \ \ \ \ \ \ \
\ \ \ \ \ \ \ \  \ \textbf{8.541}

\item[A.2]
$I_{\nu}(z)=e^{-\frac{\pi}{2}i\nu}J_{\nu}(e^{\frac{\pi}{2}i}z),
\ \ -\pi<\arg z\le\pi/2,$  \ \ \mbox{and} \\
$I_{\nu}(z)=e^{\frac{3\pi}{2}i\nu}J_{\nu}(e^{-\frac{3\pi}{2}i}z),
\ \ \ \pi/2< \arg z \le \pi.$  \ \ \ \ \ \ \ \ \ \ \ \ \ \ \ \ \ \
\ \  \ \ \ \ \ \  \ \textbf{8.406}

\item[A.3] $I_{\nu}(z)=\frac{(z/2)^{\nu}}{\Gamma(\nu+1)}
\prod_{m=1}^{\infty}\big(1+ \frac{z^2}{z_{\nu,m}^2}\big), \ \ \nu
\neq -1,-2,-3,...$,

\mbox{This is derived from the factorization A.1 with a usage of A.2.}

\item[A.3a]
$I_{\nu}(z)=(z/2)^{\nu}\sum_{k=0}^{\infty}\frac{z^{2k}}{k!\,\Gamma(\nu+k+1)}$,
\ \ \ \ \ \  \ \ \ \ \ \ \ \ \ \ \ \ \ \ \ \ \ \ \  \ \ \ \ \ \ \
\ \ \ \ \ \ \ \ \ \ \ \ \ \ \textbf{8.445}

\item[A.4]
$J_{\nu}(z)=\frac{(z/2)^{\nu}}{\Gamma(\nu+\frac{1}{2})\Gamma(\frac{1}{2})}
\int_{-1}^1 e^{izt}(1-t^2)^{\nu- \frac{1}{2}} \,dt, \ \ \Re\nu > -
\frac{1}{2},$ \ \ \ \ \ \ \ \ \ \ \ \ \ \ \ \ \ \textbf{8.41}(10)

\item[A.5]
$I_{\nu}(z)=\frac{(z/2)^{\nu}}{\Gamma(\nu+\frac{1}{2})\Gamma(\frac{1}{2})}
\int_{-1}^1 e^{\pm zt}(1-t^2)^{\nu- \frac{1}{2}}\,dt, \ \ \Re\nu >
- \frac{1}{2}$, \ \ \  \ \ \ \ \ \ \ \ \ \ \ \ \ \textbf{8.431}(1)

\item[A.6] $I_{\pm(n+ \frac{1}{2})}(z)=\frac{1}{\sqrt{2\pi\,z}}
\sum_{k=0}^n \frac{(n+k)!}{2^k k!(n-k)!z^k}
\big[(-1)^k\,e^z+(-1)^{n+1}e^{-z}\big]$, \ \ \ \ \ \textbf{8.467}

with $a_{k,n}:=(-1)^k\frac{(n+k)!}{2^k\,k!(n-k)!}$, $0\le
k\le n$, which may be rewritten as

\item[A.6a]
$I_{\pm(n+\frac{1}{2})}(z)=\sqrt{\frac{2}{\pi}}z^{-(n+\frac{1}{2})}$

$\cdot\,\big[\big(\sum_{\{k:n-k\,even\}}a_{k,n}\,z^{n-k}\big)\,\sinh z +
\big(\sum_{\{k:n-k\,odd\}}a_{k,n}\,z^{n-k}\big)\,\cosh z \big]$

\item[A.7]
$J_{n+\frac{1}{2}}(z)=(-1)^n\,z^{n+\frac{1}{2}}\sqrt{\frac{2}{\pi}}
\frac{d^n}{(z\,dz)^n} \big(\frac{\sin z}{z}\big)$ \ \ \ \ \ \ \ \
\  \ \ \ \ \ \ \ \ \ \ \ \   \ \ \ \ \ \ \ \ \ \ \textbf{8.463}(1)

[$\frac{d^n}{(z\,dz)^n}$ is the n-times composition of an operator
$ f \to \frac{1}{z}\,f'(z)$].

\item[A.7a]
$I_{n+\frac{1}{2}}(z)=z^{n+\frac{1}{2}}\sqrt{\frac{2}{\pi}}
\frac{d^n}{(z\,dz)^n}\big(\frac{\sinh z}{z}\big)$

\mbox{ It is a consequence of the formula A.7 and A.2.}

\item[A.8] $z\,\frac{d}{dz}\,I_{\nu}(z)-\nu
\,I_{\nu}(z)= z I_{\nu+1}(z)$  \ \ \  \ \ \ \ \ \ \ \ \ \ \ \ \ \ \
\ \ \ \ \ \ \ \ \ \ \ \  \ \ \ \ \ \ \ \ \ \ \ \ \textbf{8.486}(4)

\end{itemize}

\begin{center}
\textbf{REFERENCES}
\end{center}
\medskip

\noindent [1] Ph. Biane and M. Yor (1987). Variations sur une
formule de Paul L\'evy. \emph{ Ann. Inst. H. Poincar\'e Probab.
Statist.} 23(2), pp. 359-377.

\noindent [2] B. Duplantier (1989).  Areas of planar Brownian
curves, \emph{J. Phys.A: Math. Gen.} 22, 3033-3048.

\noindent [3] I. S. Gradshteyn and I. M. Ryzhik (1965).
\emph{Tables of integrals, series, and products}, Academic Press,
New York.

\noindent[4] A. Iksanov, Z. J. Jurek, B. M. Schreiber (2003). A new
factorization property of the selfdecomposale probability
measures. \emph{Ann. Prob.} (to appear)

\noindent[5] M. E. Ismail, D. H. Kelker (1979). Special functions,
Stieltjes transforms and infinite divisibility. \emph{Siam J.
Math. Anal.} 10, pp. 884-901.

\noindent [6] Z. J. Jurek (1985). Relations between the
s-selfdecomposable and selfdecomposable measures,\emph{Ann.
Probab.} 13, pp. 592-608.

\noindent [7] Z. J. Jurek (1996). Series of independent exponential
random variables. In: Proc. 7th Japan-Russia Symposium Probab.
Theory and Math. Statisitics, Tokyo 26-30 July 1995; S. Watanabe,
M. Fukushima,Yu.V. Prohorov and A.N. Shiryaev Eds. pp.174-182.
World Scientific, Singapore.

\noindent[8] Z. J. Jurek (2000). A note on gamma random variables and Dirichlet series, Stat.
Prob. Letters, 49, pp. 387-392.

\noindent [9] Z. J. Jurek (2001) Remarks on the selfdecomposability
and new examples. \emph{Demonstratio Math.}, XXXIV(2), pp. 241-250.

\noindent[10] Z. J. Jurek and J.D. Mason (1993). Operator-limit
distributions in probability theory . John Wiley and Sons, New
York.

\noindent[11] P. L\'evy (1950). Wiener's random function, and other
Laplacian random functions, \emph{Proc. 2nd Berkeley Symp. Math.
Stat. Prob.}, vol.II (Berkeley, CA: University of California
Press), 171-186.

\noindent[12] E. Lukacs (1968). Contributions to a problem of D. van
Dantzig, \emph{Theory of Probab. and Appl.} 13 (1), pp. 114-125.
(Russian Edition)

\noindent[13] B. De Meyer, B. Royenette, P. Vallois and M. Yor (2002).
On independent times and position for Brownian motion,
\emph{Revista Math. Ibero-Americana}, in print.

\noindent[14] J. Pitman and M. Yor (1981). Bessel processes and
infiniteley divisible laws. \emph{Lect. Notes in Math.} vol. 851, pp. 285-370, Springer-Verlag.

\noindent[15] D. Williams (1976). On a stopped Brownian motion formula of
H. M. Taylor, \emph{S\'eminaire de  Probab. X}, Lect. Notes in Math. 511, pp. 235-239, Springer-Verlag.

\noindent[16] M. Yor (1980). Remarques sur une formule de Paul L\'evy,
\emph{S\'eminaire de Probab. XIV}, Lect. Notes in Math. 784, pp. 343-347, Springer-Verlag.

\noindent[17] M. Yor (1989). On stochastic areas and averages of
planar Brownian motion, \emph{J. Phys. A: Math. Gen.} 22, pp. 3049-3057.

\medskip
\noindent
Author's address: Institute of Mathematics, The University of Wroc\l aw,
Pl.Grunwaldzki 2/4, 50-384 Wroc\l aw, Poland. e-mail: zjjurek@math.uni.wroc.pl

\end{document}